\documentclass{elsart}

\usepackage{amssymb,graphicx,amsmath,latexsym}
\newtheorem{theorem}{Theorem}
\newtheorem{definition}{Definition}
\newtheorem{corollary}{Corollary}
\newtheorem{lemma}{Lemma}
\def\eop{\hfill$\Box$}
\def\proof{{\em Proof: }}

\def\R{{\mathbb R}}

\def\norminf{{\|\cdot\|_{\infty}}}
\def\normeuc{{\|\cdot\|_2}}

\def\Xstar{{X^*=\{\alpha e:\alpha\in\R\}}}

\begin{document}

\begin{frontmatter}



\title{On some properties of contracting matrices}


\author{Chai Wah Wu}
\ead{chaiwahwu@member.ams.org}
\address{IBM Research Division,
Thomas J. Watson Research Center\\
P. O. Box 218,
Yorktown Heights, NY 10598, U. S. A.}

\begin{abstract}
The concepts of paracontracting, pseudocontracting and nonexpanding operators have been shown to be useful in proving convergence of asynchronous or parallel iteration algorithms.  The purpose of this paper is to give characterizations of these operators when they are linear and finite-dimensional.   
First we show that pseudocontractivity of stochastic matrices with respect to $\norminf$ is equivalent to the scrambling property, a concept first introduced in the study of inhomogeneous Markov chains. This unifies results obtained independently using different approaches.
Secondly, we generalize the concept of pseudocontractivity to set-contractivity which is a useful generalization with respect to the Euclidean norm.  In particular, we demonstrate non-Hermitian matrices that are set-contractive for $\normeuc$, but not pseudocontractive for $\normeuc$ or $\norminf$.  For constant row sum matrices we characterize set-contractivity using matrix norms and matrix graphs.  Furthermore, we prove convergence results in compositions of set-contractive operators and illustrate the differences between set-contractivity in different norms.
Finally, we give an application to the global synchronization in coupled map lattices.
\end{abstract}

\end{frontmatter}

\section{Introduction\label{sec:introduction}}
\begin{definition}[\cite{nelson:paracontractive:1987}]
Let $\|\cdot\|$ be a vector norm in ${\mathbb C}^n$.  An $n$ by $n$ matrix $B$ is {\em nonexpansive} with respect to $\|\cdot\|$ if 
\begin{equation}\label{eqn:noncontractive}
\forall x \in {\mathbb C}^n, \|Bx\|\leq \|x\|
\end{equation}
$B$ is called {\em paracontracting} with respect to $\|\cdot\|$ if 
\begin{equation}\label{eqn:paracontracting}
\forall x \in {\mathbb C}^n, Bx \neq x \Leftrightarrow \|Bx\| < \|x\|
\end{equation}
\end{definition}

It is easy to see that normal matrices with eigenvalues in the unit circle and for which $1$ is the only eigenvalue of unit norm is paracontractive
with respect to $\normeuc$.  

\begin{definition}
For a vector $x\in {\mathbb C}^n$ and a closed set $X^*$,
$y^*$ is called a {\em projection vector} of $x$ onto $X^*$ if $y^*\in X^*$ and
\[ \|x-y^*\| = \min_{y\in X^*} \|x-y\|\]
The distance of $x$ to $X^*$ is defined as $d(x,X^*) = \|x-P(x)\|$ where
$P(x)$ is a projection vector of $x$ onto $X^*$.
\end{definition}

Even though the projection vector is not necessarily unique, we write $P(x)$ when it is clear which projection vector we mean or when the choice is immaterial.
Let us denote $e = (1,\cdots , 1)^T$.
The proof of the following Lemma is relatively straightforward and thus omitted.

\begin{lemma}\label{lem:projection}
If $x\in \R^n$ and $X^* = \{\alpha e:\alpha\in \R\}$, the projection vector $P(x)$ of 
$x$ onto $X^*$ is
$\alpha e$  where:
\begin{itemize}
\item for the norm $\|\cdot\|_2$, $\alpha = \frac{1}{n}\sum_i x_i$ and $d(x,X^*) = \sqrt{\sum_i \left(x_i-\alpha\right)^2}$.
\item for the norm $\|\cdot\|_{\infty}$, $\alpha = \frac{1}{2}\left(\max_i x_i + \min_i x_i\right)$, and $d(x,X^*) = \frac{1}{2}\left(\max_i x_i - \min_i x_i\right)$.
\item for the norm $\|\cdot\|_1$, $d(x,X^*) =
\sum_{i=\lceil\frac{n}{2}\rceil+1}^n \hat{x}_i - \sum_{i=1}^{\lfloor\frac{n}{2}\rfloor}
\hat{x}_i$ and
\begin{itemize}\item for $n$ odd, $\alpha = \hat{x}_{\lceil\frac{n}{2}\rceil}$.
\item for $n$ even, $\alpha$ can be chosen
to be any number in the interval 
$[\hat{x}_{\frac{n}{2}},\hat{x}_{\frac{n}{2}+1}]$.
\end{itemize}
\end{itemize}
Here $\hat{x}_i$ are the values $x_i$ rearranged in nondecreasing order
$\hat{x}_1\leq \hat{x}_2 \leq \cdots$.
\end{lemma}

The property of paracontractivity is used to show convergence of infinite products of paracontractive matrices and this in turn is used to prove convergence in various parallel and asynchronous iteration methods \cite{bru:paracontractive:1994}.  In 
\cite{su:pseudocontractive:2001} this property is generalized to pseudocontractivity.

\begin{definition}[\cite{su:pseudocontractive:2001}]
Let $T$ be an operator on ${\mathbb R}^n$.  $T$ is {\em nonexpansive} with respect to $\|\cdot\|$ and a closed set $X^*$ if 
\begin{equation}\label{eqn:su_nonexpansive}
\forall x\in {\mathbb R}^n, x^*\in X^*, \|Tx-x^*\| \leq \|x-x^*\|
\end{equation}
$T$ is {\em pseudocontractive} with respect to $\|\cdot\|$ and $X^*$ if it is nonexpansive
with respect to $\|\cdot\|$ and $X^*$ and
\begin{equation}\label{eqn:pseudocontractive}
\forall x\not\in X^*, d(Tx,X^*) < d(x,X^*)
\end{equation}
\end{definition}
Ref. \cite{su:pseudocontractive:2001} shows that there are pseudocontractive nonnegative matrices which are not paracontractive with respect to $\norminf$ and proves a result on the convergence of infinite products of pseudocontractive matrices.   Furthermore, Ref. \cite{su:pseudocontractive:2001} studies a class of matrices for which a finite product of matrices from this class of length at least $n-1$ is pseudocontractive in $\norminf$.

The purpose of this paper is multifold.  First we show that for stochastic matrices with respect to $\norminf$ and $X^* = \{\alpha e:\alpha\in {\mathbb R}\}$, pseudocontractive matrices are equivalent to scrambling matrices and thus are simply characterized.  The concept of scrambling matrices is first introduced in the study of weak ergodicity in inhomogeneous Markov chains and
this equivalence allows us to unify several results obtained independently using these different concepts.

The second goal of this paper is to generalize pseudocontractivity by introducing the concept of set-contractivity.  We prove a convergence result of set-contractive matrices and show existence of set-contractive matrices  in $\normeuc$ that are not pseudocontractive with respect to $\normeuc$ or $\norminf$. 
We study set-contraction with respect to $\normeuc$ in terms of matrix norms and graphs of matrices.

Finally, we apply these results to the global synchronization of coupled map lattices.

We concentrate on the case where $T$ are matrices and $X^*$ is the span of the corresponding Perron eigenvector. If the Perron eigenvector is strictly positive, then as in \cite{su:pseudocontractive:2001}, a scaling operation $T\rightarrow W^{-1}TW$ where
$W$ is the diagonal matrix with the Perron eigenvector on the diagonal, transforms $T$ into a matrix for which the Perron eigenvector is $e$.  Therefore in the sequel we will focus on constant row sum matrices with
$\Xstar$.

\section{Pseudocontractivity and scrambling stochastic matrices\label{sec:scrambling}}

Scrambling matrices were first defined in \cite{hajnal:weak_ergodic:1958} to study weak ergodicity of inhomogeneous Markov chains.

\begin{definition}
A matrix $A$ is {\em scrambling} if for any pair of indices $i,j$, there exists $k$
such that $A_{ik}\neq 0$ and $A_{jk}\neq 0$.
\end{definition}

\begin{definition}
For a real matrix $A$, $\mu (A)$ is
defined as
\[ \mu(A) = \min_{j,k} \sum_{i} \min(A_{ji},A_{ki})
\]
\end{definition}

For nonnegative matrices with row sums $\leq r$, it is clear that 
$0\leq \mu(A)\leq r$ with $\mu(A) > 0$
if and only if $A$ is scrambling.

\begin{definition}
For a real matrix $A$, define $\delta(A)\geq 0$ as 
\[
\delta(A) = \max_{i,j} \sum_k \max(0,A_{ik}-A_{jk}) \geq \max_{i,j,k} (A_{ik}-A_{jk})
\]
\end{definition}

If $A$ has constant row sums, then
$\delta(A) = \frac{1}{2}\max_{i,j}\sum_k |A_{ik}-A_{jk}|$.

\begin{theorem}\label{thm:paz}
If $A$ is a matrix where each row sum is equal to or less than $r$, 
then $\delta(A) \leq r-\mu(A)$.
\end{theorem}
\proof Ref. \cite{paz:hajnal:1967} proved this for the case of stochastic matrices and the same proof applies here.\eop

\begin{theorem}\label{thm:delta}
If $A$ is a real matrix with constant row sums and $x\in \R^n$,
then $\max_i y_i - \min_i y_i \leq \delta(A)\left(\max_i x_i - \min_i x_i\right)$
where $y = Ax$.
\end{theorem}
\proof The proof is similar to the argument in \cite{paz:hajnal:1967}.
Let $x_{\max} = \max_i x_i$, $x_{\min} = \min_i x_i$, $y_{\max} = 
\max_{i} y_i$, $y_{\min} = \min_i y_i$.
\begin{equation}
\begin{array}{lcl}y_{\max}-y_{\min} & = & \max_{i,j}\sum_k \left(A_{ik}-A_{jk}\right)x_k
\\ &\leq & \max_{i,j}\left(\sum_k \max\left(0,A_{ik}-A_{jk}\right)x_{\max}
 + \sum_k \min\left(0,A_{ik}-A_{jk}\right)x_{\min}\right)
\end{array}
\end{equation}
Since $A$ has constant row sums, $\sum_k A_{ik}-A_{jk} = 0$, i.e.
\[\sum_k\max\left(0,A_{ik}-A_{jk}\right) 
+ \sum_k \min(\left(0,A_{ik}-A_{jk}\right) = 0\]
This means that 
\begin{equation}
\begin{array}{lcl}y_{\max}-y_{\min} & \leq &
\max_{i,j}\left(\sum_k\max\left(0,A_{ik}-A_{jk}\right)\right)\left(x_{\max}-x_{\min}\right)
\\ & \leq & \delta(A)\left(x_{\max}-x_{\min}\right)
\end{array}
\end{equation}
\eop

The following result shows that pseudocontractivity of stochastic matrices with respect to $\norminf$ is equivalent to the scrambling condition and thus can be easily determined.

\begin{theorem}\label{thm:scrambling-contract}
Let $A$ be a stochastic matrix.  The matrix $A$ is pseudocontractive with respect to $\norminf$ and $X^* = \{\alpha e:\alpha\in \R^n\}$ if
and only if $A$ is a scrambling matrix.
\end{theorem}
\proof
Let $x^*\in X^*$.  Then $Ax^* = x^*$ and thus
$\|Ax-x^*\|_{\infty} = \|A(x-x^*)\|_{\infty} \leq \|x-x^*\|_{\infty}$.  Thus
all stochastic matrices are nonexpansive with respect to $\norminf$ and $X^*$.
Suppose $A$ is a scrambling matrix.  Then $\mu(A) > 0$, and $\delta(A)< 1$ by Theorem \ref{thm:paz}. By Lemma \ref{lem:projection} and Theorem \ref{thm:delta}, $A$ is pseudocontractive.
Suppose $A$ is not a scrambling matrix.  Then there exists $i$,$j$ such that
for each $k$, either $A_{ik} = 0$ or $A_{jk}=0$.
Define $x$ as $x_k = 1$ if $A_{ik}> 0$ and $x_k = 0$ otherwise.  Since $A$ is stochastic,
it does not have zero rows and thus there exists $k'$ and $k''$ such that $A_{ik'} = 0$ and $A_{ik''} > 0$.
This means that $x\not\in X^*$. 
Let $y = Ax$.  Then $y_i = 1$ and $y_j = 0$.  This means that
$\max_i y_i - \min_i y_i \geq 1 = \max_i x_i - \min_i x_i$, i.e. $A$ is not pseudocontractive.\eop

With Theorem \ref{thm:scrambling-contract} several results which were shown independently can now be seen to be equivalent.  For instance, in \cite{wu:randomsynch:2006} it was shown that for stochastic matrices with positive diagonal entries and whose interaction digraph\footnote{The directed graph of a square matrix $A$ is defined as the graph with an edge from vertex $i$ to vertex $j$ if and only if $A_{ij} \neq 0$.  The interaction digraph of a matrix $A$ is obtained from the directed graph of $A$ by reversing the orientation of all the edges, i.e. it is the graph of $A^T$.} contains
a spanning directed tree a finite product of $n-1$ or more such matrices is scrambling. 
In \cite{wu:reducible:2005} it was shown that such matrices are irreducible or 1-reducible\footnote{A matrix is 1-reducible if after simultaneous row and column permutation it can be written in the form
$\left(\begin{array}{cccc}B_{11} & B_{12} & \cdots & \\ & B_{22} & B_{23} & \cdots \\ & & \ddots & \\ & & &B_{kk}\end{array}\right)$ such that $B_{ii}$ are irreducible and for each $i < k$, there exists $j> i$ such that $B_{ij}\neq 0$.} and this result in \cite{wu:randomsynch:2006}
then mirrors Proposition 3.3 in \cite{su:pseudocontractive:2001}.

In \cite{lubachevsky:async:1986} the convergence of a class of asynchronous iteration algorithms was shown by appealing to results about scrambling matrices.  In \cite{su:pseudocontractive:2001} this result is proved using the framework of pseudocontractions.  Theorem \ref{thm:scrambling-contract} shows that these two approaches are essentially the same.

\section{Set-nonexpansive and set-contractive operators}
Consider the stochastic matrix 
\[A = \left(\begin{array}{ccc}0.5& 0 & 0.5\\0.5&0.5&0\\0.5&0&0.5\end{array}\right)\]
The matrix $A$ is not
pseudocontractive with respect to the Euclidean norm $\normeuc$ and $X^* = \{\alpha e:\alpha\in\R\}$ since $\|A\|_2 = 1.088 > 1$.  On the other hand, $A$ satisfies Eq. (\ref{eqn:pseudocontractive})\footnote{This can be shown using Theorem \ref{thm:ev}.}.  This motivates us to define the following generalization of pseudocontractivity:

\begin{definition}\label{def:x-expand}
Let $X^*$ be a closed set in $\R^n$.  An operator $T$ on $\R^n$ is {\em set-nonexpansive} with respect to $\|\cdot\|$ and $X^*$ if
\[ \forall x\in \R^n, d(Tx,X^*) \leq d(x,X^*) \]
An operator $T$ on $\R^n$ is {\em set-contractive} with respect to $\|\cdot\|$ and $X^*$
if it is set-nonexpansive with respect to $\|\cdot\|$ and $X^*$ and
\[ \forall x\not\in X^*, d(Tx,X^*) < d(x,X^*).
\]
The {\em set-contractivity} of an operator $T$ is defined as
\[ c(T) = \sup_{x\not\in X^*}\frac{d(Tx,X^*)}{d(x,X^*)} \geq 0\]
\end{definition}

There is a dynamical interpretation to Definition \ref{def:x-expand}.  If we consider the operator $T$ as a discrete-time dynamical system, then $T$ being set-nonexpansive and set-contractive
imply that $X^*$ is a globally nonrepelling invariant set and a globally attracting set of the dynamical system respectively \cite{wiggins_tam_book_90}.

\begin{lemma}\label{lem:one}
$T$ is set-nonexpansive with respect to $\|\cdot\|$ and $X^*$ if and only if $T(X^*)\subseteq X^*$ and $c(T)\leq 1$. 
If $T$ is set-contractive with respect to $\|\cdot\|$ and $X^*$, then the fixed points of $T$ is a subset of $X^*$.
If $T_1(X^*)\subseteq X^*$, then $c(T_1\circ T_2) \leq c(T_1)c(T_2)$.
\end{lemma}
\proof The first statement is true by definition. The proof of the second statement is the same as in Proposition 2.1 in \cite{su:pseudocontractive:2001}.
Suppose $T_1(X^*)\subseteq X^*$.  Let $x\not\in X^*$.  If $T_2(x)\in X^*$, then 
$d(T_1\circ T_2(x),X^*) = 0$.  If $T_2(x)\not\in X^*$, then 
$d(T_1\circ T_2(x)) \leq c(T_1)d(T_2(x),X^*) \leq c(T_1)c(T_2)d(x,X^*)$.\eop

\begin{lemma}\label{lem:ct}
Let $X^*$ be a closed set such that $\alpha X^*\subseteq X^*$ for all $\alpha\in\R$. If $T$ is linear and $T(X^*)\subseteq X^*$,
then $c(T) = \sup_{\|x\|=1,P(x) = 0} d(T(x),X^*)$.  
\end{lemma}
\proof Let $\epsilon = \sup_{\|x\|=1,P(x) = 0} d(T(x),X^*)$.  Clearly $\epsilon \leq c(T)$.  For $x\not\in X^*$, $0$ is a projection vector of $x-P(x)$.
Since $T(P(x))\in X^*$, this implies that $d(T(x),X^*) = d(T(x-P(x)),X^*) \leq \epsilon \|x-P(x)\| = \epsilon d(x,X^*)$, i.e. $\epsilon \geq c(T)$.
\eop

\begin{lemma}\label{lem:contractivity-linear}
Let $X^*$ be a closed set such that $\alpha X^*\subseteq X^*$ for all $\alpha\in\R$.  An set-nonexpansive matrix $T$ is set-contractive with respect to $X^*$ if and only if $c(T) < 1$.
\end{lemma}
\proof One direction is clear.  Suppose $T$ is set-contractive.  By compactness
\[ \sup_{\|x\|=1,P(x) = 0} d(T(x),X^*) = \epsilon < 1 \]
and the conclusion follows from Lemma \ref{lem:one} and Lemma \ref{lem:ct}.
\eop

If $T$ is nonexpansive with respect to $\|\cdot \|$ and $X^*$, then
\[ \|Tx-P(Tx)\| \leq \|Tx-P(x)\| \leq \|x-P(x)\| \]
and $T$ is set-nonexpansive.  Thus set-contractivity is more general than pseudocontractivity.
However, they are equivalent for stochastic matrices with respect to $\norminf$ and $X^*=\{\alpha e:\alpha\in\R\}$.

\begin{lemma}
With respect to  $\norminf$ and $X^*=\{\alpha e:\alpha\in\R\}$,
a stochastic matrix $T$ is pseudocontractive if and only if it is set-contractive.
\end{lemma}
\proof Follows from the fact that a stochastic matrix is nonexpansive with respect to $\norminf$ and $X^*=\{\alpha e:\alpha\in\R\}$.\eop

\begin{definition}[\cite{horn-johnson:matrix_analysis:1985}]
A vector norm $\|\cdot\|$ on $\R^n$ is {\em monotone} if
\[\left\| \left( x_1, \cdots , x_n\right)^T \right\| \leq
\left\| \left( y_1, \cdots , y_n\right)^T  \right\|
\]
for all $x_i$ and $y_i$ such that $|x_i|\leq |y_i|$.
A vector norm  $\|\cdot\|$ on $\R^n$ is {\em weakly monotone} if
\[\left\| \left( x_1, \cdots, x_{k-1} , 0, x_{k+1}, \cdots , x_n\right)^T \right\| \leq
\left\| \left( x_1, \cdots, x_{k-1} , x_k, x_{k+1}, \cdots , x_n\right)^T  \right\|
\]
for all $x_i$ and $k$.
\end{definition}

The next result gives a necessary condition of set-contractivity of a matrix in terms of its graph.

\begin{theorem}\label{thm:sdt}
Let $A$ be a constant row sum matrix with row sums $r$ such that $|r|\geq 1$.  If
$A$ is set-contractive with respect to a weakly monotone vector norm $\|\cdot \|$ 
and $\Xstar$, then the interaction digraph of $A$ contains a spanning directed tree.
\end{theorem}
\proof 
If the interaction digraph $A$ does not have a spanning directed tree, it was shown in \cite{wu:reducible:2005} that after simultaneous row and column permutation, $A$ can be written as a block upper triangular matrix:
\[ A = \left(\begin{array}{cccc} * & * & * & * \\ & \ddots & * & * \\ & & A_1& 0 \\
& & & A_2\end{array}\right) \]
where $*$ are arbitrary entries and $A_1$ and $A_2$ are $m_1$ by $m_1$ and $m_2$ by $m_2$ square irreducible matrices respectively.  Define $x = (0,\dots , 0,-a_1 e_1, a_2 e_2)^T \not\in X^*$,
where $e_1$ and $e_2$ are vectors of all $1$'s of length $m_1$ and $m_2$ respectively.
Let $z = (0,\dots , 0, e_3)^T$ where $e_3$ is the vector of all $1$'s of length $m_1+m_2$
and $Z^* = \{\alpha z: \alpha \in \R\}$.
Note that the set of projection vectors of a fixed vector $x$ to $Z^*$ is a convex connected set.
Let $\alpha z$ be a projection vector of $x$ to $Z^*$.  Suppose that for $a_1 = a_2 \neq 0$, $\alpha \neq 0$.  Since $-\alpha z$ is a projection vector of $-x$ 
to $Z^*$ and $\alpha$ (or at least a choice of $\alpha$) depends continuously on $a_1$ and $a_2$, by varying $a_1$ to $-a_1$ and varying $a_2$ to $-a_2$, $\alpha$ changes to $-\alpha$.  This means that we can find $a_1$ and $a_2$ not both zero, such that $0$ is a projection vector of $x$ to $Z^*$.  In this case $x\not\in X^*$ and by weak monotonicity 
$d(x,Z^*) = d(x,X^*) = \|x\|$.
It is clear that $y = Ax$ can be written as 
\[y = \left(\begin{array}{c}*\\\vdots \\ *\\  -ra_1e_1\\ ra_2e_2
\end{array}\right)\]
Let $\beta e$ be a projection vector of $y$ onto $X^*$. 
By the weak monotonicity of the norm,
\[ d(y,X^*) = \|y-\beta e\| \geq
\left(\begin{array}{c} 0\\\vdots \\ 0\\  (-ra_1-\beta)e_1\\ (ra_2-\beta)e_2
\end{array}\right)
= r\left(x-\frac{\beta}{r}z\right)
\] 
Since $0$ is a projection vector of $x$ onto $Z^*$
\[ d(y,X^*) \geq |r| d(x,Z^*) \geq d(x,X^*) \]
Thus $A$ is not set-contractive.\eop

\subsection{max-norm}

\begin{theorem}\label{thm:contractinf}
Let $A$ be a matrix with constant row sum $r$. Then $c(A) = r-\mu(A)$ with respect to $\norminf$ and $\Xstar$.  In particular,
the matrix $A$ is set-nonexpanding with respect to $\norminf$ and $\Xstar$
if and only if $r-\mu(A) \leq 1$.
The matrix
$A$ is set-contractive with respect to $\norminf$ and $\Xstar$ if and only if $r-\mu(A) < 1$.
\end{theorem}
\proof
$c(A)\leq r-\mu(A)$ follows from Lemma \ref{lem:projection}, Theorem \ref{thm:paz} and Theorem \ref{thm:delta}.  Since $c(A)\geq 0$, $c(A) = r - \mu(A)$ if $r-\mu(A) = 0$.
Therefore we assume that $r-\mu(A) > 0$.
Let $j$ and $k$ be such that $\mu(A) = \sum_i \min\left(A_{ji},A_{ki}\right)$.
Define $x$ such that $x_i = 1$ if $A_{ji} < A_{ki}$ and $x_i = 0$ otherwise.
Since $r - \mu(A) > 0$, $x$ is not all $0$'s or all $1$'s, i.e. $x\not\in X^*$.
Let $y = Ax$.  Then by Lemma \ref{lem:projection}
\[\begin{array}{lcl} 2d(y,X^*) \geq  y_k-y_i & = & \sum_{i,A_{ji} < A_{ki}} A_{ki}-A_{ji} \\
&=&
\sum_i A_{ki} - \sum_{i,A_{ji} \geq A_{ki}} A_{ki} - \sum_{i,A_{ji} < A_{ki}} A_{ji} 
\\ &
=& r - \mu(A) \end{array}
\]
Since $2d(x,X^*) = 1$, it follows that $c(A) \geq r-\mu(A)$.
\eop

\subsection{Euclidean norm}

The following result characterizes set-contractivity of matrices with respect to $\normeuc$ in terms of 
matrix norms.

\begin{theorem} \label{thm:ev}
Let $A$ be an $n$ by $n$ constant row sum matrix and $K$ be an $n$ by $n-1$ matrix whose columns form a orthonormal basis of $e^{\bot}$.  
Then $c(A) = \left\|AK\right\|_2$ with respect to $\normeuc$ and $\Xstar$.
In particular $\left\|AK\right\|_2 \leq 1$
if and only if $A$ is set-nonexpanding with respect to $\normeuc$ and  $X^*=\{\alpha e:\alpha\in\R\}$.  Similarly,  $\left\|AK\right\|_2 < 1$ if and only if 
$A$ is set-contracting with respect to $\normeuc$ and  $X^*=\{\alpha e:\alpha\in\R\}$.
\end{theorem}
\proof  Define $J = ee^T$ as the $n$ by $n$ matrix of all $1's$. Note that $\|x\|_2 = \|Kx\|_2$ and $JK = 0$.  Let $B = A-\frac{1}{n}J$.
Then
\[ \|AK\|_2 = \|BK\|_2 =
\max_{\|x\|_2=1}\|BKx\|_2 = \max_{\|Kx\|_2=1}\|BKx\|_2 =
\max_{x\bot e,\|x\|_2=1}\|Bx\|_2
\]
By Lemma \ref{lem:projection} $P(x) = \frac{1}{n}Jx$ and $d(Ax,X^*) = \|Bx\|$.  Since $A$ has constant row sums, $A(X^*)\subseteq X^*$ and by Lemma \ref{lem:ct}
$c(A) = \max_{P(x) = 0,\|x\|_2=1} d(Ax,X^*) = \max_{P(x) = 0,\|x\|_2=1} \left\|Bx\right\|_2$.
Since $P(x) = 0$ if and only if $x\bot e$, this means that
$c(A) = \left\|AK\right\|_2$.
\eop

\subsection{weighted Euclidean norm}
\begin{definition}
Given a positive vector $w$, the weighted $2$-norm $\|\cdot\|_w$ is defined as
 \[ \|x\|_w = \sqrt{\sum_{i} w_i x_i^2} \]
\end{definition}

\begin{theorem} \label{thm:ev2}
Let $A$ be an $n$ by $n$ constant row sum matrix and $K$ be as defined in Theorem \ref{thm:ev}. Let $w$ be a positive vector such that $\max_i w_i = 1$ and $W = \mbox{diag}(w)$.
Then $c(A) \leq \left\|W^{\frac{1}{2}}AW^{-1}K\right\|_2$ with respect to $\|\cdot\|_w$ and $\Xstar$.
\end{theorem}
\proof The proof is similar to Theorem \ref{thm:ev}.
Define $J_w = \frac{ew^T}{\sum_i w_i}$ and $B = A-J_w$.
Note that $J_wW^{-1}K = 0$.
Then 
\[ \begin{array}{lcl} \|W^{\frac{1}{2}}AW^{-1}K\|_2 &= &\|W^{\frac{1}{2}}BW^{-1}K\|_2 
\\ &=&  
 \max_{\|Kx\|_2=1}\|W^{\frac{1}{2}}BW^{-1}Kx\|_2 \\ &=&
\max_{x\bot e,\|x\|_2=1}\|W^{\frac{1}{2}}BW^{-1}x\|_2\end{array}\]
Now $x\bot e$ if and only if $W^{-1}x \bot w$.  Since 
$\|x\|_2 = \|W^{-\frac{1}{2}}x\|_w$, this means that
$\|W^{\frac{1}{2}}AW^{-1}K\|_2 = \max_{x\bot w,\|W^{\frac{1}{2}}x\|_w = 1}\|W^{\frac{1}{2}}Bx\|_2$.
Since $\max_i w_i = 1$, this means that $\|W^{\frac{1}{2}}x\|_w 
= \sqrt{\sum_i (w_ix_i)^2} \leq \|x\|_w$ and thus
\[\|W^{\frac{1}{2}}AW^{-1}K\|_2 \geq \max_{x\bot w,\|x\|_w = 1}\|W^{\frac{1}{2}}Bx\|_2
\]

It is straightforward to show that $P(x) = J_wx$ and thus
 $d(Ax,X^*) = \|Bx\|_w = \|W^{\frac{1}{2}}Bx\|_2$.  Since $A$ has constant row sums, $A(X^*)\subseteq X^*$ and by Lemma \ref{lem:ct}
$c(A) = \max_{P(x) = 0,\|x\|_w=1} d(Ax,X^*) = \max_{P(x) = 0,\|x\|_w=1} \left\|W^{\frac{1}{2}}Bx\right\|_2$.
Since $P(x) = 0$ if and only if $x\bot w$, this means that
$c(A) \leq \left\|W^{\frac{1}{2}}AW^{-1}K\right\|_2$.\eop

Note that the matrix $A$ in Theorem \ref{thm:contractinf}, Theorem \ref{thm:ev} and Theorem \ref{thm:ev2} is not necessarily nonnegative or stochastic.

\subsection{examples}
The matrix
\[A_1 = \left(\begin{array}{ccc}1.1& 0.0 & 0.0\\0.6&0.5&0\\0.6&0&0.5\end{array}\right)\]
is set-contracting with respect to $\norminf$ and $\Xstar$ since $\mu(A_1) = 0.6$
and $c(A_1) = 1.1-\mu(A_1) = 0.5 < 1$.  It is not pseudocontracting with respect to
$\norminf$ and $\Xstar$ since $\|A_1\|_{\infty} = 1.1 > 1$.

The stochastic matrix
\[A_2 = \left(\begin{array}{ccc}0.4& 0.3 & 0.3\\0&1&0\\0&0&1\end{array}\right)\]
is set-nonexpanding with respect to $\normeuc$ and $\Xstar$ since $\left\|A_2K\right\|_2 = 1$ but it is not nonexpanding with respect to $\normeuc$ and $X^*=\{\alpha e:\alpha\in\R\}$ since $\|A_2\|_2 > 1$.  Furthermore, 
Theorem \ref{thm:sdt} shows that $A_2$ is not set-contractive with respect to any
weakly monotone norm and $X^*$.

The stochastic matrix
\[A_3 =\left(\begin{array}{ccc}1& 0 & 0\\0.5&0.5&0\\0&0.5&0.5\end{array}\right)\] is set-contractive with respect to $\normeuc$ and $\Xstar$ since $\left\|A_3K\right\|_2 = 0.939$.  Since $\|A_3\|_2 > 1$ it is not nonexpanding nor pseudocontractive with respect to $\normeuc$ and $X^*$. It is also not pseudocontractive 
with respect to $\norminf$ and $X^*$ since it is not scrambling.

The stochastic matrix
\[A_4 =\left(\begin{array}{ccc}1& 0 & 0\\0.9&0.1&0\\0.1&0.1&0.8\end{array}\right)\]
has an interaction digraph that contains a spanning directed tree.
However, it is not set-nonexpanding with respect to  $\normeuc$ and $\Xstar$
since $\left\|A_4K\right\|_2 = 1.125 > 1$.  
This shows that
the converse of Theorem \ref{thm:sdt} is not true for $\normeuc$.\footnote{Theorem
\ref{thm:scrambling-contract} shows that the converse of Theorem \ref{thm:sdt} is false as well for stochastic matrices with respect to $\norminf$ and $X^*$.}  
On the other hand, $A_4$ is set-contractive with respect to $\norminf$ and $X^*$ since $A_4$ is a scrambling matrix.  Furthermore, $A_4$ is set-contractive with respect to
$\|\cdot\|_w$ and $X^*$ for $w = (1,0.2265,1)^T$ since $\|W^{\frac{1}{2}}A_4W^{-1}K\|_2 < 1$.

Next we show some convergence results for dynamical systems of the form
$x(k+1) = T_kx(k)$ where some $T_k$'s are 
set-contractive operators.

\begin{theorem}\label{thm:conv}
Let $\{T_k\}$ be a sequence of set-nonexpansive operators with respect to $\|\cdot\|$ 
and $X^*$ and suppose that
\[ \lim_{k\rightarrow\infty} \prod_k c(T_k) = 0\]
Let $x(k+1) = T_kx(k)$. For any initial
vector $x(0)$, $\lim_{k\rightarrow\infty} d(x(k),X^*) = 0$.
\end{theorem}
\proof From Lemma \ref{lem:one}, 
$c(\prod_k T_k) \leq \prod_k c(T_k) \rightarrow 0$ as $k\rightarrow\infty$ and the conclusion follows.\eop

\begin{theorem}
Let $\Xstar$ and $\{A_k\}$ be a sequence of $n$ by $n$ constant row sum nonnegative matrices
such that
\begin{itemize}
\item the diagonal elements are positive;
\item all nonzero elements are equal to or larger than $\epsilon$;
\item the row sum is equal to or less than $r$.
\end{itemize}
If 
$r^{n-1}-\epsilon^{n-1} < 1$ and for each $k$, the interaction digraph of $A_k$ contains a spanning directed tree, then
$\lim_{k\rightarrow \infty} d(x(k),X^*) = 0$ where $x(k+1) = A_kx(k)$.
\end{theorem}
\proof As discussed above, products of $n-1$ matrices $A_k$
is scrambling.  By definition, since each $A_k$ has nonzero elements equal to or larger than $\epsilon$,  the nonzero elements of this product, denoted as $P$, will be equal to or larger than $\epsilon^{n-1}$.  This means that $\mu(P) \geq \epsilon^{n-1}$ and thus
$\delta(P) \leq r^{n-1}-\epsilon^{n-1} < 1$ since $P$ has row sums $\leq r^{n-1}$.
Therefore $P$ is set-contractive with respect to 
$\norminf$ and $X^*$ with $c(P) \leq r^{n-1}-\epsilon^{n-1} < 1$. The result then follows from Theorem \ref{thm:conv}.\eop

The following result shows existence of linear operators $B_k$ and vectors $x_k^*\in X^*$ such that $x(k+1) = B_kx(k) + x^*_k$ has the same dynamics as
$x(k+1) = T_kx(k)$.  In particular, for $y(k+1) = B_ky(k)$ and $x(k+1) = T_kx(k)$,
$d(y(k),X^*) = d(x(k),X^*)$ for all $k$.

\begin{theorem}
$T$ is a set-nonexpansive operator with respect to $\norminf$ and $\Xstar$ 
if and only if for each $x\in\R^n$ there exists a stochastic matrix $B$ and a vector
$x^*\in X^*$such that 
$T(x) = Bx + x^*$.

$T$ is a set-contractive operator with respect to $\norminf$ and $\Xstar$ 
if and only if for each $x\in\R^n$ there exists a scrambling stochastic matrix $B$ and 
a vector $x^*\in X^*$ such that 
$T(x) = Bx + x^*$.
\end{theorem}

\proof One direction of both statements follows from Theorem \ref{thm:scrambling-contract}.  Suppose $T$ is set-nonexpansive and fix $x\in \R^n$.
Define $x^* = P(T(x)) - P(x)$ which is a vector in $X^*$.
Let $y = T(x)-x^*$. Then $P(y) = P(T(x))-x^* = P(x)$ and by Lemma \ref{lem:projection}, 
\[ \min_i x_i \leq \min_i y_i \leq \max_i y_i \leq \max_i x_i \]
and thus there exists a stochastic matrix $B$ such that $Bx = y$.

If $T$ is set-contractive, then for $x\in X^*$, we can choose $B = \frac{1}{n}ee^T$
and
$T(x)-Bx\in X^*$.
For $x\not\in X^*$, $d(x,X^*)< d(T(x),X^*)$.  Define $x^*$ and $y$ as before and  
we see that
\[ \min_i x_i < \min_i y_i \leq \max_i y_i < \max_i x_i \]
If $x_{i'} = \min_i x$, then it is clear that we can pick $B$ with $Bx = y$ such that the $i'$-th column of $B$ is positive, i.e. $B$ is scrambling.\eop

It can be beneficial to consider set-contractivity with respect to different norms.
For instance, consider $x(k+1) = A_k x(k)$ where $A_k$ are matrices that are not pseudocontractive with respect to $\norminf$ and $\Xstar$ and whose diagonal elements are $0$.  Since the diagonal elements are not positive, the techniques in \cite{su:pseudocontractive:2001} cannot be used to show that products of $A_k$ are pseudocontractive with respect to $\norminf$ and $\Xstar$.  However, it is possible that $A_k$ are set-contractive with respect to
a different norm and thus convergence of $x(k)$ can be obtained by studying set-contractivity using this norm.
For instance, the stochastic matrix
\[ A = \left(\begin{array}{ccc} 0 & 0.5 & 0.5\\ 1 & 0 & 0\\ 0.5 & 0.5 & 0\end{array}\right)\]
has zeros on the diagonal and is not pseudocontractive with respect to $\norminf$ and $\Xstar$ since $A$ is not scrambling.  On the other hand, $A$
is set-contractive with respect to $\normeuc$ and $\Xstar$ since
$\left\|AK\right\|_2 = 0.939 < 1$.

For a set of constant row sum matrices $A_k$ and $x(k+1) = A_kx(k)$, a lower bound for the exponential rate at which $x(k)$ approach $\Xstar$ is $-\mbox{ln}(c(A))$.
The above examples show that there are matrices for which this rate is $0$ for $\norminf$ and positive for $\normeuc$ and other matrices for which 
the rate is positive and $0$ for $\norminf$ and $\normeuc$ respectively.

On the other hand, even though set-contractivity depends on the norm used, the equivalence of norms on $\R^n$ and Lemma \ref{lem:contractivity-linear} provides the following result.

\begin{theorem}
Let $X^*$ be a closed set such that $\alpha X^*\subseteq X^*$ for all $\alpha\in\R$.
and let $H$ be a compact set of set-contractive matrices with respect to $\|\cdot\|_p$ and $X^*$.  Then there exists $m$ such that a product of $m$
matrices in $H$ is set-contractive with respect to $\|\cdot\|_q$.
\end{theorem}

\begin{corollary}
Let $H$ be a compact set of stochastic set-contractive matrices with respect to $\|\cdot\|_p$ and $\Xstar$.  Then a sufficiently long product of matrices in $H$ is scrambling.
\end{corollary}

\section{Weak ergodicity of inhomogeneous Markov chains}
In Section \ref{sec:scrambling} we noted the connection between set-contractivity with respect to $\norminf$ and weak ergodicity in inhomogeneous Markov chains.
In this section we elaborate on this connection.
A sequence of stochastic matrices $A_i$ is {\em weakly ergodic}
if for each $r$, $\delta\left(A_rA_{r+1}\cdots A_{r+k}\right)\rightarrow 0$ as 
$k\rightarrow\infty$.

In \cite{seneta:markov:1973} a {\em coefficient of ergodicity} is defined as a continuous function $\mu$ on the set of $n$ by $n$ stochastic matrices such that $0\leq \mu(A)\leq 1$.
A coefficient of ergodicity $\mu$ is {\em proper} if
\[ \mu(A) = 1 \Leftrightarrow A = ev^T \quad \mbox{for some probability vector $v$.}\]

Seneta \cite{seneta:markov:1973} gives the following necessary and sufficient conditions for weak ergodicity generalizing the arguments by Hajnal.

\begin{theorem}\label{thm:seneta}
Suppose $\mu_1$ and $\mu_2$ are coefficients of ergodicity such that $\mu_1$ is proper
and the following equation is satisfied for some constant $C$ and all $k$,
\begin{equation}\label{eqn:mult} 1-\mu_1(S_1S_2\cdots S_k) \leq C\prod_{i=1}^k (1-\mu_2(S_i)) \end{equation}
where $S_i$ are stochastic matrices.
Then a sequence of stochastic matrices $A_i$ is weakly ergodic if there exists
a strictly increasing subsequence $\{i_j\}$ such that
\begin{equation}\label{eqn:ergodic} \sum_{j=1}^{\infty} \mu_2 (A_{i_j+1}\cdots A_{i_{j+1}}) = \infty\end{equation}
Conversely, if $A_i$ is a weakly ergodic sequence, and $\mu_1$, $\mu_2$ are both proper
coefficients of ergodicity
satisfying Eq. (\ref{eqn:mult}),
then Eq. (\ref{eqn:ergodic}) is satisfied for some strictly increasing sequence $\{i_j\}$.
\end{theorem}

Define $H$ as the set of stochastic matrices that are set-nonexpansive with respect to 
a norm $\|\cdot\|$ and $\Xstar$.  For $\norminf$, $H$ is the set of stochastic matrices.
Let us define $\mu_c(A) = 1-c(A)$.  Then $\mu_c$ is a proper coefficient of ergodicity when restricted to $H$.  This can be seen as follows.  Clearly
$0\leq \mu_c(A) \leq 1$.  If $A = ev^T$, then $Ax\in X^*$ and thus $c(A) = 0$ and 
$\mu_c(A) = 1$.  If $A \neq ev^T$, then there exists $i,j,k$ such that $A_{ik}\neq A_{jk}$.
Let $x$ be the $k$-th unit basis vector.  Then $(Ax)_i\neq (Ax)_j$, i.e. $d(Ax,X^*) > 0$,
$c(A) > 0$ and $\mu_c(A) < 1$.
By choosing $\mu_1 = \mu_2 = \mu_c$, Eq. (\ref{eqn:mult}) is satisfied
with $C=1$ by Lemma \ref{lem:one}.
Thus we have shown that a sufficient and necessary condition for a sequence
of matrices in $H$ to be weakly ergodic is 
\[\sum_{j=1}^{\infty} 1- c(A_{i_j+1}\cdots A_{i_{j+1}}) = \infty\]
for some strictly increasing subsequence $\{i_j\}$.

\section{Application to the synchronization of coupled map lattices}
Coupled map lattices \cite{kaneko:cml_overview_1992} have been studied extensively and have been shown to exhibit complex behavior \cite{kaneko:global_coupled_1989,manneville:random_coupling_1992}.  Recently, synchronization in coupled map lattice has attracted considerable attention \cite{belykh:1dlattice:1996,gade:cml_synch_1996,wu:cml-1998,jost:cml:2002,lu:synch-discrete:2004}.  We show here how set-contractivity can be useful in studying synchronization in coupled map lattices.

Given a map $f_k:\R\rightarrow\R$, consider state variables $x_i\in \R$ which evolve according to
$f_k$ at time $k$: $x_i(k+1) = f_k(x(k))$.  By coupling the output of these maps we obtain a coupled map lattice where each state evolves as:
\[ x_i(k+1) = \sum_{j} a_{ij}(k) f_k(x_j(k)) \]
This can be rewritten as
\begin{equation}\label{eqn:cml} x(k+1) = A_k F_k(x(k))\end{equation}
where $x(k) = (x_1(k),\dots, x_n(k))^T \in \R^n$ and
$F_k(x(k)) = (f_k(x_1(k)),\dots , f_k(x_n(k)))^T$.  We assume that $A_k$ is a constant row sum matrix for all $k$. The map $f_k$
depends on $k$, i.e. we allow the map in the lattice to be time varying.  Furthermore, we do not require $A_k$ to be a nonnegative matrix.  We say the coupled map lattice in Eq. (\ref{eqn:cml}) {\em synchronizes}
if $\lim_{k\rightarrow\infty} |x_i(k)-x_j(k)| = 0$ for all $i$ and $j$, i.e.
$x(k)$ approaches the synchronization manifold $\Xstar$ 
as $k\rightarrow\infty$.
If the row sum of $A_k$ is $1$, then this means that at synchronization,
each state $x_i$ in the lattice 
exhibits dynamics of the uncoupled map $f_k$, i.e. if $x(h)\in X^*$, then for all $k\geq h$,
$x(k)\in X^*$ and $x_i(k+1) = f_k(x_i(k))$.

We are now ready to state our synchronization result:
\begin{theorem} \label{thm:cml}
Let $\rho_k$ be the Lipschitz constant of $f_k$.
If $\lim_{k\rightarrow\infty} \prod_k c(A_k)\rho_k = 0$, where $c(A_k)$ is the set-contractivity with respect to $\Xstar$ and a monotone norm, then the coupled map lattice in Eq. (\ref{eqn:cml})
synchronizes.
\end{theorem}
\proof 
\[ \|F_k(x(k))- P(F_k(x(k))\| \leq 
\|F_k(x(k)) - F_k(P(x(k))\| \leq
\rho_k \|x(k)-P(x(k)\|\]
where the last inequality follows from monotonicity of the norm.  
This implies that $c(F_k)\leq \rho_k$
and the result follows from Theorem \ref{thm:conv}.\eop

Thus we can synchronize the coupled map lattice if we can find matrices $A_k$ and a norm
such that the contractivities $c(A_k)$ are small enough.

\begin{corollary}
Let $\rho_k$ be the Lipschitz constant of $f_k$. If $\sup_k r(A_K)-\mu(A_K)-\frac{1}{\rho_k} < 0$, 
then Eq. (\ref{eqn:cml}) synchronizes\footnote{Here $r(A)$ denotes the row sum of the matrix $A$.}.
\end{corollary}
\proof Follows by applying Theorem \ref{thm:cml} to set-contractivity with respect to 
$\norminf$.\eop

\bibliography{chua_ckt2,chaos,secure,consensus,synch,misc,stability,cml,algebraic_graph,graph_theory,control,optimization,adaptive,top_conjugacy,ckt_theory,cnn2,matrices,chaos_comm,markov}
\end{document}